\documentclass[final,1p,times]{elsarticle}
\usepackage{amssymb}
\usepackage{amsmath}
\usepackage{amsthm}

\newtheorem{theorem}{Theorem}[section]

\numberwithin{equation}{section}

\newtheorem{lemma}[theorem]{Lemma}

\theoremstyle{definition}

\journal{www.arxiv.org}

\begin{document}

\begin{frontmatter}



\title{A note on boundary differentiability of solutions of nondivergece elliptic equations with unbounded drift\tnoteref{t1}}

%% \address[label2]{<address>} and 11101324
\tnotetext[t1]{The author was supported by NSFC 11401460 and CSC 201506285016.}
\author[rvt]{Yongpan Huang\corref{cor1}}
\ead{huangyongpan@gmail.com}
\cortext[cor1]{Corresponding author}
\address[rvt]{School of Mathematics and Statistics, Xi'an Jiaotong University, Xi'an 710049, China}

\begin{abstract}
Boundary differentiability is shown for solutions of nondivergence elliptic equations with unbounded drift.
\end{abstract}

\begin{keyword}
Boundary differentiability, Elliptic equations
\end{keyword}

\end{frontmatter}

% \linenumbers

\section{Introduction}

In this paper, we will study the boundary differentiability of strong solutions of elliptic equation with unbounded lower order coefficients. Suppose that $u\in W^{2,n}_{loc}(Q^+_1)\cap C(\overline{Q^+_1})$ satisfies
\begin{eqnarray}\label{EQ1}
 \left\{\begin{array}{l} Lu:=-a_{ij}(x)D_{ij}u+b_i(x)D_iu=f(x) \;\; in  \;\;\ Q^+_1;\\
 u(x)=0  \;\; on  \;\;\  T_1.
\end{array}
\right.
\end{eqnarray}
We use the summation convention over repeated indices and the
notations $D_i:= \frac{\partial}{\partial x_i}$; $D_{ij} := D_iD_j$. We assume that $a_{ij}$, $b_i$ and $f$ are measurable
functions on $Q^+_1$, $\mathbf{b}=(b_1,b_2,...,b_n)$,  the matrix $(a_{ij}(x))_{n\times n}$ is symmetric and satisfies the uniformly elliptic condition
\begin{equation}\label{EQ2}
\lambda|\xi|^2\leq  a_{ij}(x)\xi_i\xi_j\leq \lambda^{-1}|\xi|^2,\ \ \forall \xi\in \mathbb{R}^n,\ \ a.e. x\in \ Q^+_1,
\end{equation}
with a constant $\lambda\in(0,1]$, and $b_i, f\in L^n(Q^+_1)$. Throughout the paper, we denote $W(\Omega):=W^{2,n}_{loc}(\Omega)\cap C(\overline{\Omega})$ while $\Omega$ is a bounded domain in $\mathbb{R}^n$.

 As for the boundary regularity of nondivergence elliptic equations: If the drift term $|\mathbf{b}|$ is bounded, Krylov \cite{KRY1} showed that the solution is $C^{1,\alpha}$ along the boundary if $\partial\Omega$ is $C^{1,1}$; Lieberman  \cite{GML} gave a more general estimates; Wang \cite{WANG} proved a similar pointwise result as in \cite{KRY1} by an iteration method that will be adopted in this paper; Li and Wang in \cite{LW1,LW2} showed the boundary differentiability of solutions of elliptic equations on convex domains. If $|\mathbf{b}|$ is unbounded, Ladyzhenskaya and Ural'tseva \cite{LU} proved boundary $C^{1,\alpha}$ estimate of elliptic and parabolic inequalities on $W^{2,q}$ domain with $\mathbf{b}\in L^q$, $\Phi\in L^q$, $q>n$ and nonlinear term $\mu_1|Du|^2$; Safonov \cite{MVS2} obtained the the Hopf-Oleinik lemma for elliptic equations and gave the counterexample which indicated that the Dini condition on $b_n$ can not be removed for our theorem; Nazarov  \cite{AIN} proved the Hopf-Oleinik Lemma and boundary gradient estimate under minimal restrictions on lower-order coefficients; In \cite{HQZ} the boundary differentiability is shown for strong solution of nondivergence elliptic equation  $|\mathbf{b}|$ and $f$  satisfying Dini condition. Since the Hopf Oleinik Lemma and boundary Lipschitz Estimate \cite{AIN} hold for solution of (1.1) only need $b_n$ satisfies the Dini condition, it is natural conjecture that whether the boundary differentiability of solutions at $0$ is true while $b_n$  satisfying Dini condition at $0$. In the following, we will show that the result is correct. Some related results concerning Dini continuity can be found in \cite{BUR,KOV1,KOV2,MVS1}.\\

The following Alexandroff-Bakelman-Pucci maximum principle and Harnack inequality are our main tools.
\begin{theorem}(\cite{GT,MVS2})
Let $\Omega$ be a bounded domain in $\mathbb{R}^n$, and let $u$ be a function in $W(\Omega)$ such that $Lu \leq f$ in $\Omega$. Suppose that the matrix $(a_{ij}(x))_{n\times n}$ is symmetric and satisfies the uniformly elliptic condition (\ref{EQ2}), and $b_i, f\in L^{n}(\Omega)$. Then
\begin{equation}\label{EQ3}
\sup_{\Omega}u\leq \sup_{\partial\Omega}u+N diam\Omega\cdot e^{N||\mathbf{b}||_{L^n(\Omega)}^n}||f^{+}||_{L^n(\Omega)},
\end{equation}
where
\begin{equation}\label{EQ4}
||\mathbf{b}||_{L^n(\Omega)}=\Big(\int_{\Omega}|\mathbf{b}|^ndx\Big)^{\frac{1}{n}},\ \ \ \mathbf{b}=(b_1,b_2,...,b_n),
\end{equation}
and $N$ is a positive constant depending only on $n$ and $\lambda$.
\end{theorem}

\begin{theorem}\label{thm1.2}(Harnack Inequality)
Let $u$ be a nonnegative function in $W(B_{8})$, $Lu=f$ in $B_{8}$  and $b_i, f\in L^n(B_{8})$. There exists a positive constant $\epsilon_0$ depending only on $\lambda$ and $n$, such that if $||\mathbf{b}||_{L^n(B_{8})}\leq\epsilon_0$, then
\begin{equation}\label{EQ5}
\sup\limits_{B_1}u\leq C(\inf\limits_{B_1}u+||f||_{L^n(B_{8})}),
\end{equation}
where $C$ is constant depending only on $\lambda$ and $n$.
\end{theorem}
Theorem \ref{thm1.2} follows from the the proof in \cite{MVS2} clearly. The most important thing is that the quantity $||\mathbf{b}||_{L^n}$ is scaling invariant(see Remark 1.4 in \cite{MVS2}) and the Harnack constant is invariant in the iteration procedure.

\begin{flushleft}
\textbf{Notations.}
\end{flushleft}
\noindent
$\{\vec{e}_i\}_{i=1}^n,$\ \ \ \ \ \ \ \ \ \ \ \ \ \ \ \ \ \ \ \ \  \ \  \ \  the standard basis of $\mathbb{R}^n$.\\
$|x|~:=\sqrt{\sum\limits_{i=1}^{n}x_i^2},$\ \ \ \ \ \ \ \ \ \ \ the Euclidean norm of $x=(x_1,x_2,...,x_n)\in \mathbb{R}^n$.\\
$B_r:=\{x\in \mathbb{R}^{n}: \ |x|<r\}$.\\
$T_r:=\{x'\in \mathbb{R}^{n-1}: \ |x'|<r\}.$\\
$Q_r:=T_r\times(-r,r).$\\
$Q_r^+:=T_r\times(0,r).$\\
$||f||_{L^n(\Omega)}:=\Big(\int_{\Omega}|f(x)|^ndx\Big)^{\frac{1}{n}}.$\\
$W(\Omega):=W^{2,n}_{loc}(\Omega)\cap C(\overline{\Omega}).$\\

\begin{theorem}\label{thm1.5}
Assume that
(1)$u\in W(Q^+_1)$, $u|_{T_1}=0$, $Lu=f$ in $Q^+_1$;
(2)$f\in L^n(Q^+_1)$  and $\int_0^1\frac{||f||_{L^n(Q^+_r)}}{r}dr<\infty$;
(3)$\mathbf{b}\in L^n(Q^+_1)$ and $\int_0^1\frac{||b_n||_{L^n(Q^+_r)}}{r}dr<\infty$.
Then $u$ is differentiable at $0$.
\end{theorem}

\section{Proof of Theorem}
By standard normalization, it is enough for us to prove the following Theorem \ref{thm2.1} instead of proving Theorem \ref{thm1.5}.
\begin{theorem}\label{thm2.1}
Assume that\\
(1)$u\in W(Q^+_1)$, $||u||_{L^{\infty}(Q^+_1)}\leq 1$,  $Lu=f$ in $Q^+_1$,\ \ and\  \  $u|_{T_1}=0$;\\
(2)$f\in L^n(Q^+_1)$ with $||f||_{L^n(Q^+_1)}\leq1$ and $\int_0^1\frac{||f||_{L^n(Q^+_r)}}{r}dr\leq1$;\\
(3)$\mathbf{b}$ and $b_n$ satisfy
\begin{equation}\label{Eq2.1}
||\mathbf{b}||_{L^n(Q^+_1)}\leq\min\{\frac{\delta}{4A_2(\frac{4M}{\mu}+1)},1,\epsilon_0\}\ \ \ and\ \ \ \int_{0}^{1}\frac{||b_n||_{L^n(Q^+_r)}}{r}dr\leq\min\{1, \frac{\delta\ln\frac{1}{\delta}}{16MA_3}\},
\end{equation}
where $\epsilon_0$ is the constant in Theorem \ref{thm1.2}, and $\delta$, $M$, $\mu$, $A_2$ and $A_3$ are constants in Lemma \ref{lem2.2}.\\
Then $u$ is differentiable at $0$.
\end{theorem}

\begin{lemma}\label{lem2.2}
There exist positive constant $\delta(<1)$,$\mu(<1)$, $M$, $A_1$, $A_2$ and $A_3$ depending only on $\lambda$ and $n$.
If
\begin{equation}
kx_n-B\leq u(x)\leq Kx_n+B \ \ \mbox{in}\ \ Q^+_1,
\end{equation}
for some constants $k$, $K$ and $B(\geq0)$ with $k\leq K$, then there
exist constants $\tilde{k}$ and $\tilde{K}$ such that for $x\in Q^+_{\delta}$,
\begin{eqnarray}
&&\tilde{k}x_n-
A_1||f||_{L^{n}(Q^+_1)}-A_2(K-k+B)||\mathbf{b}||_{L^n(Q^+_1)})-A_3(|K|+|k|)||b_n||_{L^n(Q^+_1)}\nonumber\\
&&\leq u(x)\leq \tilde{K}x_n+A_1||f||_{L^{n}(Q^+_1)}+A_2(K-k+B)||\mathbf{b}||_{L^n(Q^+_1)})+A_3(|K|+|k|)||b_n||_{L^n(Q^+_1)},
\end{eqnarray}
where either
\begin{equation}
\tilde{k}=k-2MB+\mu(K-k) \ \ \mbox{and}\ \ \tilde{K}=K+2MB,
\end{equation}
or
\begin{equation}
\tilde{k}=k-2MB \ \ \mbox{and}\ \ \tilde{K}=K+2MB-\mu(K-k).
\end{equation}
Obviously, we have $\tilde{k}\leq\tilde{K}$.
\end{lemma}

\proof[Proof of Lemma 2.2]~
We prove the following claim first.

\begin{flushleft}
\textbf{Claim.}  There exist positive constants $M$, $\tilde{\delta}$ and $C_1$ depending only on $\lambda$ and $n$, such
that for any $x\in Q^+_{\tilde{\delta}}$,
\begin{equation}
\begin{aligned}
&\quad(k-2MB)x_n-C_1(|k|||b_n||_{L^n(Q^+_1)}+B||\mathbf{b}||_{L^n(Q^+_1)}+||f||_{L^n(Q^+_1)})\\
&\leq u(x)\leq(K+2MB)x_n+C_1(|K|||b_n||_{L^n(Q^+_1)}+B||\mathbf{b}||_{L^n(Q^+_1)}+||f||_{L^n(Q^+_1)}).
\end{aligned}
\end{equation}
\end{flushleft}
\proof~~Let $M\geq\sqrt{n-1}(1+\frac{2\sqrt{n-1}}{\lambda})$ and $\epsilon(>0)$ be
small enough, such that
\begin{equation}
4-(1+\epsilon)(2+\epsilon)(M-1)^{\epsilon}\geq0.
\end{equation}

Let
\begin{equation}
\tilde{\delta}=\frac{1}{M}(\leq\frac{1}{3\sqrt{n-1}}),\ \ \delta=\frac{\tilde{\delta}}{2M}
\end{equation}
and
\begin{equation}
\tilde{\psi}(x)=2\Big(\frac{x_n}{\tilde{\delta}}\Big)-\Big(\frac{x_n}{\tilde{\delta}}\Big)^2+\frac{\lambda^2}{2(n-1)}\sum_{i=1}^{n-1}\left(\left(\frac{|x_i|}{\tilde{\delta}}-1\right)^{+}\right)^{2+\epsilon}.
\end{equation}
The barrier function $\tilde{\psi}(x)$ is $C^{2}$ and satisfies
the following conditions:
\begin{eqnarray}
 \left\{\begin{array}{l}
(1)\tilde{\psi}(x)\geq1\ \mbox{on}\ \{x\in \mathbb{R}^n:\ |x'|\leq1, x_n=\tilde{\delta}\};\\[3pt]
(2)\tilde{\psi}(x)\geq0\ \mbox{on}\ \{x\in \mathbb{R}^n:\ |x'|\leq1, x_n=0\};\\[3pt]
(3)\tilde{\psi}(x)\geq1\ \mbox{on}\ \{x\in \mathbb{R}^n:\ |x'|=1, 0\leq x_n\leq\tilde{\delta}\};\\[3pt]
(4)-a^{ij}(x)D_{ij}\tilde{\psi}(x)\geq0\ \ \mbox{a.e.} \ \mbox{in}\ \ \{x\in \mathbb{R}^n:|x'|< 1,\ \ \ 0< x_n<\tilde{\delta}\};\\[3pt]
(5)\tilde{\psi}(x)\leq\frac{2x_n}{\tilde{\delta}}\ \ \ \mbox{in}\
\ Q^+_{\tilde{\delta}}.
\end{array}
\right.
\end{eqnarray}

It follows that
\begin{eqnarray}
 \left\{\begin{array}{l}
L\big(kx_n-B\tilde{\psi}(x)-u(x)\big)\leq b_iD_i(kx_n-B\tilde{\psi}(x))-f(x) \;\; in  \;\;\ \tilde{Q};\\[10pt]
kx_n-B\tilde{\psi}(x)-u(x)\leq 0 \;\; on
\;\;\partial\tilde{Q},
\end{array}
\right.
\end{eqnarray}
where
$\tilde{Q}=T_1\times(0,\tilde{\delta}).$

According to the Alexandroff-Bakelman-Pucci maximum principle, we have
\begin{equation}
kx_n-B\tilde{\psi}(x)-u(x)\leq C_1(|k|||b_n||_{L^n(Q^+_1)}+B||\mathbf{b}||_{L^n(Q^+_1)}+||f||_{L^n(Q^+_1)}))\ \ \ \mbox{in}\ \ \ \tilde{Q},
\end{equation}
where $C_1$ is a constant depending only on $\lambda$ and $n$.

By (2.10)(5)(i.e.
$\tilde{\psi}(x)\leq\frac{2x_n}{\tilde{\delta}}=2Mx_n\ \ \mbox{in}\ Q^+_{\tilde{\delta}}$), we have
\begin{equation}
(k-2MB)x_n-C_1(|k|||b_n||_{L^n(Q^+_1)}+B||\mathbf{b}||_{L^n(Q^+_1)}+||f||_{L^n(Q^+_1)})\leq u(x)\ \ \ \mbox{in}\ \ \ Q^+_{\tilde{\delta}}.
\end{equation}
As in (2.11), we also have
\begin{eqnarray}
 \left\{\begin{array}{l}
L\big(u(x)-Kx_n-B\tilde{\psi}(x)\big)\leq -b_iD_i(Kx_n+B\tilde{\psi(x)})+f(x)\;\; in  \;\;\ \tilde{Q};  \\[10pt]
u(x)-Kx_n-B\tilde{\psi}(x)\leq 0 \;\; on
\;\;\partial\tilde{Q}.
\end{array}
\right.
\end{eqnarray}

According to the Alexandroff-Bakelman-Pucci maximum principle, we have
\begin{equation}
u(x)-Kx_n-B\tilde{\psi}(x)\leq \ C_1(|K|||b_n||_{L^n(Q^+_1)}+B||\mathbf{b}||_{L^n(Q^+_1)}+||f||_{L^n(Q^+_1)})\ \ \ \mbox{in}\ \ \tilde{Q},
\end{equation}
where $C_1$ is a constant depending only on $\lambda$ and $n$ .
Combining (2.15) and (2.10)(5), we get
\begin{equation}
u(x)\leq (K+2MB)x_n+C_1(|K|||b_n||_{L^n(Q^+_1)}+B||\mathbf{b}||_{L^n(Q^+_1)}+||f||_{L^n(Q^+_1)})\ \ \
\mbox{in}\ \ Q_{\tilde{\delta}}.
\end{equation}
By (2.13) and (2.16), the claim follows clearly.\qed\\

Let $\Gamma=\{\delta
\vec{e}_n+T_{M\delta}\}$. Next, we will show (2.3) according to two cases:   $u(\delta
e_n)\geq\frac{1}{2}(K+k)\delta$ and
$u(\delta e_n)<\frac{1}{2}(K+k)\delta$,
corresponding to which (2.4) and (2.5) will hold respectively. Since the proofs of these two cases are similar,  we will only show the proof for the case: $u(\delta
e_n)\geq\frac{1}{2}(K+k)\delta$.\\

Let $v(x)=u(x)-(k-2MB)x_n+C_1(|k|||b_n||_{L^n(Q^+_1)}+B||\mathbf{b}||_{L^n(Q^+_1)}+||f||_{L^n(Q^+_1)})$.
Then
\begin{equation}
v(\delta
e_n)\geq\big(\frac{K-k}{2}+2MB\big)\delta+C_1(|k|||b_n||_{L^n(Q^+_1)}+B||\mathbf{b}||_{L^n(Q^+_1)}+||f||_{L^n(Q^+_1)}).
\end{equation}
Since $v(x)\geq0$ for $x\in Q^+_{\tilde{\delta}}$, from (2.17) and the Harnack
inequality, it follows that
\begin{equation}
\sup_{\Gamma}v(x)\leq C_2\big(\inf_{\Gamma}
v(x)+|k|||b_n||_{L^{n}(Q^+_1)}+B||b_n||_{L^{n}(Q^+_1)}+||f||_{L^n(Q^+_1)}\big),
\end{equation}
where $C_2(\gg1)$ is a constant depending only on $\lambda$ and $n$. Combining (2.17),(2.18) and $v(x)\geq 0$, we
have
\begin{equation}
\inf_{\Gamma}v(x)\geq\Big\{\frac{1}{C_2}\Big((\frac{K-k}{2}+2MB)\delta\Big)+(\frac{C_1}{C_2}-1)(|k|||b_n||_{L^n(Q^+_1)}+B||\mathbf{b}||_{L^n(Q^+_1)}+||f||_{L^n(Q^+_1)}\Big\}^+:=
a.
\end{equation}

Let
\begin{eqnarray}
\psi(x)=\frac{1}{2}\left(\Big(\frac{x_n}{\delta}\Big)
+\Big(\frac{x_n}{\delta}\Big)^2\right)-\frac{\lambda^2}{4(n-1)}\sum_{i=1}^{n-1}\left(\left(\frac{|x_i|}{\delta}-1\right)^{+}\right)^{2+\epsilon},
\end{eqnarray}
where $\epsilon$ satisfies (2.7).

The barrier function $\psi(x)$ is $C^{2}$ and satisfies the
following conditions :
\begin{eqnarray}
 \left\{\begin{array}{l}
(1)\psi(x)\leq1\ \mbox{on}\ \{x\in \mathbb{R}^n:\ |x'|\leq M\delta, x_n=\delta\};\\[3pt]
(2)\psi(x)\leq0\ \mbox{on}\ \{x\in \mathbb{R}^n:\ |x'|\leq M\delta, x_n=0\};\\[3pt]
(3)\psi(x)\leq0\ \mbox{on}\ \{x\in \mathbb{R}^n:\ |x'|= M\delta, 0< x_n<\delta\};\\[3pt]
(4)-a_{ij}(x)D_{ij}\psi(x)\leq0\ \ \mbox{a.e.\ in}\ \  \{x\in \mathbb{R}^n:|x'|< M\delta,\ 0<x_n<\delta\}; \\[3pt]
(5)\psi(x)\geq\frac{x_n}{2\delta}\ \ \ \mbox{in}\ \ Q^+_\delta;\\[3pt]
%(6)\psi(x)\leq\frac{x_n}{\delta}\ \ \mbox{in}\ \{x\in \mathbb{R}^n:\ |x'|\leq
%M\delta, 0\leq x_n\leq\delta\}.
\end{array}
\right.
\end{eqnarray}
It follows that
\begin{eqnarray}
 \left\{\begin{array}{l}
L\big(a\psi(x)-v(x)\big)\leq b_iD_i(a\psi(x)+(k-2MB)x_n)-f(x) \;\;\;\; in  \;\;\ \tilde{\tilde{Q}};\\[10pt]
a\psi(x)-v(x)\leq 0\;\; on
\;\;\partial\tilde{\tilde{Q}},
\end{array}
\right.
\end{eqnarray}
where
$\tilde{\tilde{Q}}=T_{M\delta}\times(0,\delta).$

According to the Alexandroff-Bakelman-Pucci maximum principle,
\begin{equation}
a\psi(x)-v(x)\leq C_3(K-k+B)||\mathbf{b}||_{L^n(Q^+_1)}+C_4|k|||b_n||_{L^n(Q^+_1)}+C_5||f||_{L^n(Q^+_1)}\ \ \ \mbox{in}\ \ \ \tilde{\tilde{Q}},
\end{equation}
where $C_3,C_4,C_5$ are constants depending only on $\lambda$ and $n$, and we have used $K-k\geq0$.

From (2.21)(5), it follows that for each $x\in Q^+_\delta$,
\begin{align}
a\psi(x)&\geq\frac{a}{2\delta}x_n\nonumber\\
&\geq\frac{\frac{(K-k)\delta}{2C_2}-|k|||b_n||_{L^n(Q^+_1)}-B||\mathbf{b}||_{L^n(Q^+_1)}-||f||_{L^n(Q^+_1)}}{2\delta}x_n\nonumber\\
&\geq\frac{K-k}{4C_2}x_n-|k|||b_n||_{L^n(Q^+_1)}-B||\mathbf{b}||_{L^n(Q^+_1)}-||f||_{L^n(Q^+_1)}.
\end{align}
Combining (2.23) and (2.24), we have that  for  each $x\in Q^+_\delta$,
\begin{equation}
\begin{aligned}
u(x)&\geq a\psi(x)+(k-2MB)x_n-(C_1+C_4)(|k|||b_n||_{L^n(Q^+_1)}+B||\mathbf{b}||_{L^n(Q^+_1)})\\
&\quad~ -C_3(K-k+B)||\mathbf{b}||_{L^n(Q^+_1)}-(C_1+C_5) ||f||_{L^n(Q^+_1)}\\
&\geq\Big(k-2MB+\frac{1}{4C_2}(K-k)\Big)x_n-(C_1+C_4+1)|k|||\mathbf{b}||_{L^n(Q^+_1)}.\\
&\quad~ -(C_1+C_3+C_4+1)(K-k+B)||\mathbf{b}||_{L^n(Q^+_1)}-(C_1+C_5+1) ||f||_{L^n(Q^+_1)}.
\end{aligned}
\end{equation}
Let
\begin{equation}
\mu=\frac{1}{4C_2},\ \ \ A_1=C_1+C_5+1,\ \ \ A_2=C_1+C_3+C_4+1\ \ and\ \ A_3=C_1+C_4+1.
\end{equation}
Combining (2.16),(2.25) and (2.26) , we have that (2.3) and (2.4) hold.\qed\\

By induction, the following lemma is a direct consequence of Lemma 2.2.
\begin{lemma}
There exist sequences
$\{k_m\}_{m=0}^{\infty}$ and $\{K_m\}_{m=0}^{\infty}$, and nonnegative sequence $\{B_m\}_{m=0}^{\infty}$ with $k_0=K_0=0$ , $B_0=1$, and for
$m=0,1,2,...,$
\begin{align}
B_{m+1}&=A_1\delta^m||f||_{L^{n}(Q^+_{\delta^m})}+A_2\delta^m(K_m-k_m+\frac{B_m}{\delta^m})||\mathbf{b}||_{L^n(Q^+_{\delta^m})}+A_3\delta^m(|K_m|+|k_m|)||b_n||_{L^n(Q^+_{\delta^m})}\nonumber,
\end{align}
and
$$k_{m+1}=k_m-2M\frac{B_m}{\delta^m}+\mu(K_m-k_m)\ \ \mbox{and}\ \ K_{m+1}=K_m+2M\frac{B_m}{\delta^m},$$
or
$$k_{m+1}=k_m-2M\frac{B_m}{\delta^m}\ \ \mbox{and}\ \ K_{m+1}=K_m+2M\frac{B_m}{\delta^m}-\mu(K_m-k_m),$$
such that
\begin{equation}
k_mx_n-B_m\leq u(x)\leq K_mx_n+B_m\ \ \mbox{in}\ \ Q^+_{\delta^m},
\end{equation}
where $\delta$, $\mu$, $M$, $A_1$ and $A_2$ are positive constants given by Lemma 2.2.
\end{lemma}

%\begin{remark}
%It seems that there are no hopes for us to prove the convergence of $\sum\limits_{i=0}^{\infty}\frac{B_i}{\delta^i}$ because of the term $||\mathbf{b}||_{L^n(Q^+_{\delta^m})}$ appears, but the estimate for partial sum of $\{K_i-k_i\}_{i=1}^{\infty}$ helps us to prove it. We have no explicit estimate as in \cite{MW1}.
%\end{remark}

Now we present the proof of Theorem 2.1.

\proof[Proof of Theorem 2.1]~Let
$\{B_m\}_{m=0}^{\infty}$, $\{k_m\}_{m=0}^{\infty}$ and
$\{K_m\}_{m=0}^{\infty}$ be defined by Lemma 2.3. We will show the proof by the following three claims.

\begin{flushleft}
\textbf{Claim 1.} $\sum\limits_{m=0}^{\infty}\frac{B_m}{\delta^m}$ is convergent.
\end{flushleft}

\proof~~Firstly, notice that we take $K_0=k_0=0$ and $B_0=1$, then by induction, we have $K_m\geq k_m$ for all $m\geq0$.

For $m\geq0$, we define $S_m=\sum\limits_{i=0}^{m}\frac{B_i}{\delta^i}.$  For any $m\geq0$, since $$ K_{m+1}\leq K_m+2M\frac{B_m}{\delta^m}\ \ \mbox{and}\ \ K_0=0,$$
we have
\begin{equation}
K_{m+1}\leq 2MS_m\ \ \mbox{for any}\ \ \ m\geq0.\nonumber
\end{equation}
Similarly, we have
\begin{equation}
k_{m+1}\geq -2MS_m\ \ \mbox{for any}\ \ \ m\geq0.\nonumber
\end{equation}
Hence,
\begin{equation}
|K_{m+1}|+|k_{m+1}|\leq 4MS_m\ \ \mbox{for any}\ \ \ m\geq0.
\end{equation}
Now we consider the term $K_m-k_m$.
By Lemma 2.3, for any $m\geq0$,
\begin{equation}
K_{m+1}-k_{m+1}\leq (1-\mu)(K_{m}-k_{m})+4M\frac{B_m}{\delta^m}.\nonumber
\end{equation}
Since $K_0=k_0=0$, by iteration, we have that for any $m\geq0$,
\begin{equation}
K_{m+1}-k_{m+1}\leq \sum\limits_{i=0}^m\frac{4MB_i}{\delta^i}(1-\mu)^{m-i}.
\end{equation}
It follows that for $m\geq1$,
\begin{equation}
\sum\limits_{j=0}^m(K_{j}-k_{j})\leq \sum\limits_{j=1}^m\sum\limits_{i=0}^{j-1}\frac{4MB_i}{\delta^i}(1-\mu)^{j-1-i}=\sum\limits_{j=0}^{m-1}\sum\limits_{i=0}^{j}\frac{4MB_i}{\delta^i}(1-\mu)^{j-i}.\nonumber
\end{equation}
By changing the order of summation, we have
\begin{equation}
\sum\limits_{j=0}^{m-1}\sum\limits_{i=0}^{j}\frac{4MB_i}{\delta^i}(1-\mu)^{j-i}=\sum\limits_{i=0}^{m-1}\sum\limits_{j=i}^{m-1}\frac{4MB_i}{\delta^i}(1-\mu)^{j-i}.\nonumber
\end{equation}
By
\begin{equation}
\sum\limits_{j=i}^{\infty}(1-\mu)^{j-i}=\frac{1}{\mu},\nonumber
\end{equation}
we have that for $m\geq1$,
\begin{equation}
\sum\limits_{j=0}^m(K_{j}-k_{j})\leq \frac{4M}{\mu}\sum\limits_{i=0}^{m-1}\frac{B_i}{\delta^i}=\frac{4M}{\mu}S_{m-1}.
\end{equation}

Since
$$\frac{B_{i+1}}{\delta^{i+1}}=\frac{A_1}{\delta}||f||_{L^{n}(Q^+_{\delta^i})}+\frac{A_2}{\delta}(K_i-k_i+\frac{B_i}{\delta^i})||\mathbf{b}||_{L^n(Q^+_{\delta^i})}+\frac{A_3}{\delta}(|K_i|+|k_i|)||b_n||_{L^n(Q^+_{\delta^i})},$$
for any $i\geq1$, combining the above identity with (2.28), we obtain
\begin{equation}
\frac{B_{i+1}}{\delta^{i+1}}\leq \frac{A_1}{\delta}||f||_{L^{n}(Q^+_{\delta^i})}+\frac{A_2}{\delta}(K_i-k_i+\frac{B_i}{\delta^i})||\mathbf{b}||_{L^n(Q^+_{\delta^i})}
+\frac{4MA_3}{\delta}S_{i-1}||b_n||_{L^n(Q^+_{\delta^i})}.
\end{equation}

It follows from (2.30) and (2.31) that for any $m\geq1$,
\begin{align}
\sum\limits_{i=1}^m\frac{B_{i+1}}{\delta^{i+1}}&\leq \frac{A_1}{\delta}\sum\limits_{i=1}^m||f||_{L^{n}(Q^+_{\delta^i})}+\frac{A_2}{\delta}\sum\limits_{i=1}^m(K_i-k_i+\frac{B_i}{\delta^i})||\mathbf{b}||_{L^n(Q^+_{\delta^i})}
+\frac{4MA_3}{\delta}\sum\limits_{i=1}^mS_{i-1}||b_n||_{L^n(Q^+_{\delta^i})}\nonumber\\
&\leq \frac{A_1}{\delta}\sum\limits_{i=1}^m||f||_{L^{n}(Q^+_{\delta^i})}+\frac{A_2}{\delta}(\frac{4M}{\mu}+1)S_{m+1}||\mathbf{b}||_{L^n(Q^+_{1})}
+\frac{4MA_3}{\delta}S_{m+1}\sum\limits_{i=1}^m||b_n||_{L^n(Q^+_{\delta^i})}.\nonumber
\end{align}

Since
\begin{align}
\sum_{i=1}^{\infty}\frac{4MA_3}{\delta}||b_n||_{L^n(Q^+_{\delta^i})}&\leq\frac{4MA_3}{\delta \ln{\frac{1}{\delta}}}\int_{0}^{1}\frac{||b_n||_{L^n(Q^+_r)}}{r}dr\leq\frac{1}{4},\nonumber\\
\frac{A_2}{\delta}(\frac{4M}{\mu}+1)||\mathbf{b}||_{L^n(Q^+_{1})}&\leq \frac{1}{4},\nonumber
\end{align}
and
\begin{equation}
\sum_{i=1}^{\infty}||f||_{L^{n}(Q^+_{\delta^{i}})}\leq \frac{1}{\ln{\frac{1}{\delta}}}\int_{0}^{1}\frac{||f||_{L^n(Q^+_r)}}{r}dr,\nonumber
\end{equation}
it follows that
\begin{align}
S_{m+1}-S_{1}&=\sum_{i=1}^{m}\frac{B_{i+1}}{\delta^{i+1}}\leq\frac{A_1}{\delta\ln{\frac{1}{\delta}}}\int_{0}^{1}\frac{||f||_{L^n(Q^+_r)}}{r}dr +\frac{1}{2}S_{m+1}.\nonumber
\end{align}

Therefore for all $m\geq1$,
\begin{equation}
S_{m+1}\leq\frac{2A_1}{\delta\ln{\frac{1}{\delta}}}\int_{0}^{1}\frac{||f||_{L^n(Q^+_r)}}{r}dr+2S_{1}\leq\frac{4A_1}{\delta\ln{\frac{1}{\delta}}}+2S_{1}\leq\frac{4A_1}{\delta\ln{\frac{1}{\delta}}}+\frac{2(A_1+A_2)}{\delta}+2,\nonumber
\end{equation}
where we used $||f||_{L^{n}(Q^+_{1})}\leq1$, $\int_{0}^{1}\frac{||f||_{L^{n}(Q^+_{r})}}{r}dr\leq1$ and $||\mathbf{b}||_{L^{n}(Q^+_{1})}\leq1$.
Then $\{S_m\}_{m=0}^{\infty}$ is a uniformly bounded sequence. It follows that $\sum\limits_{m=0}^{\infty}\frac{B_m}{\delta^m}$ is convergent.  This completes the proof of Claim 1.\qed\\

\begin{flushleft}
\textbf{Claim 2.} $$\lim\limits_{m\rightarrow\infty}K_m=\lim\limits_{m\rightarrow\infty}k_m=\theta.$$
\end{flushleft}

\proof~~It follows from Claim 1 that $\{K_m\}_{m=0}^{\infty}$ and $\{k_m\}_{m=0}^{\infty}$ are uniformly bounded.

Since
$$K_{m+1}-K_m\leq 2M\frac{B_m}{\delta^m}=2MS_m-2MS_{m-1}\ \ \mbox{for}\ \ m\geq1, $$
we obtain
\begin{equation}
K_{m+1}-2MS_m\leq K_m-2MS_{m-1}\ \ \mbox{for}\ \ m\geq1.\nonumber
\end{equation}
It follows that $\{K_m-2MS_{m-1}\}_{m=1}^{\infty}$ is a bounded nonincreasing sequence and $\lim\limits_{m\rightarrow\infty}(K_m-2MS_{m-1})$ exists. Hence $\lim\limits_{m\rightarrow\infty}K_m$ exists. Let $\lim\limits_{m\rightarrow\infty}K_m=\theta$.

Since
\begin{equation}
\sum\limits_{j=0}^m(K_{j}-k_{j})\leq \frac{2M}{\mu}\sum\limits_{i=0}^{\infty}\frac{B_i}{\delta^i}<+\infty,\ \ \ \ \ \forall\ \ m\geq1,\nonumber
\end{equation}
we have  $\sum\limits_{j=0}^\infty(K_{j}-k_{j})$ is convergent. It follows that $\lim\limits_{m\rightarrow\infty}(K_m-k_m)=0$.  Hence  $$\lim\limits_{m\rightarrow\infty}K_m=\lim\limits_{m\rightarrow\infty}k_m=\theta$$

This completes the proof of Claim 2.\qed\\

\begin{flushleft}
\textbf{Claim 3.} Let $\theta$ be given by Claim 2. Then for each
$m=0,1,2,...$, there exist $C_m$ such that
$\lim\limits_{m\rightarrow\infty}C_m=0$ and that $|u(x)-\theta x_n|\leq
C_m\delta^m$ for any $x\in Q^+_{\delta^m}$.
\end{flushleft}

\proof~~For any $m\geq0$ and any $x\in Q^+_{\delta^m}$, we have
\begin{equation}
|u(x)-\theta x_n|\leq(|K_m-\theta|+|k_m-\theta|)|x_n|+\frac{B_m}{\delta^m}\leq(|K_m-\theta|+|k_m-\theta|+\frac{B_m}{\delta^m}\})\delta^m.\nonumber
\end{equation}

Let
$C_m=|K_m-\theta|+|k_m-\theta|+\frac{B_m}{\delta^m}$. It follows that for any $m\geq0$ and any $x\in Q^+_{\delta^m}$,
\begin{equation}
|u(x)-\theta x_n|\leq C_m\delta^m,\nonumber
\end{equation}
and
\begin{equation}
\lim\limits_{m\rightarrow\infty}C_m=0.\nonumber
\end{equation}

This completes the proof of Claim 3.\qed\\

By Claim 3, we deduce that $u(x)$ is differentiable at $0$ with derivative $\theta \vec{e}_n$. This completes the proof of Theorem 2.1. \qed


\begin{thebibliography}{00}
\bibitem{AN1}D.E.Apushkinskaya, A.I.Nazarov; \textit{Boundary estimate for the first-order derivatives of a solution to a nondivergent parabolic equation with composite right-hand side and coefficients of lower-order derivatives}, Journal of Mathematical Sciences, \textbf{77}(4) (1995), 3257-3276.
\bibitem{AN2}D.E.Apushkinskaya, A. I. Nazarov; \textit{A counterexample to the Hopf-Oleinik lemma (elliptic case)}, Analysis $\&$ PDE, \textbf{9}(2) (2016), 439--458.
\bibitem{BUR}C. Burch; \textit{The Dini condition and regularity of weak solutions of elliptic equations}, Journal of Differential Equations, \textbf{30}(1970), 308--323.
\bibitem{GT}D.Gilbarg, N.S.Trudinger; \textit{Elliptic partial differential equations of second order},  2nd
ed., Springer-Verlag, Berlin, (1983).
\bibitem{LU}O.A.Ladyzhenskaya, N.N. Ural'tseva; \textit{Estimates on the boundary of a domain for the first derivatives of functions satisfying an elliptic or parabolic inequality}, Trudy Mat. Inst.Steklov, \textbf{179} (1988), 102--125 (in Russian). English transl. in Proc. Steklov Inst. Math., \textbf{179} (1989), 109--135.
\bibitem{KOV1}J. Kovats; \textit{Fully nonlinear elliptic equations and the Dini condition}, Comm. Partial Differential Equations, \textbf{22} (1997), 1911--1927.
\bibitem{KOV2}J. Kovats; \textit{Dini--Campanato spaces and applications to nonlinear elliptic equations}, Electron. J. Differential Equations, \textbf{37} (1999), 1--20.
\bibitem{KRY1}N.V.Krylov; \textit{Boundedly inhomogeneous elliptic and parabolic equations in a domain}, Izvestia Akad. Nauk. SSSR, \textbf{47} (1983), 75--108 (in Russian). English translation in Math. USSR Izv., \textbf{22} (1984), 67--97.
\bibitem{GML}G.M.Lieberman; \textit{The Dirichlet problem for quasilinear elliptic equations with continuously differentiable boundary data}, Comm. in Partial Differential Equations, \textbf{11} (1986), 167--229.
\bibitem{LW1} D.Li, L.Wang; \textit{Boundary differentiability of solutions of elliptic equations on convex domains}, Manuscripta Mathematica, \textbf{121} (2006), 137--156.
\bibitem{LW2}D.Li, L.Wang; \textit{Elliptic equations on convex domains with nonhomogeneous Dirichlet boundary conditions}, Journal of Differential Equations, \textbf{246} (2009), 1723--1743.
%\bibitem{MW1}F.Ma, L.Wang; \textit{Boundary first order derivatives estimates for fully nonlinear elliptic equations}, Journal of Differential Equations, \textbf{252} (2012), 988--1002.
\bibitem{HQZ}Y. Huang, Q. Zhai, S. Zhou; \textit{Boundary regularity for strong solution of nondivergence elliptic equations with unbounded drift},   Electronic Journal of Differential Equations,  Vol. 2019 (2019), No. 39, 1-16.
\bibitem{AIN}A.I.Nazaro; \textit{A centennial of Zaremba-Hopf-Oleinik Lemma, SIAM Jouranl of mathematcial analysis}, SIAM J. Math. Anal., \textbf{44} (2012), 437--453.
\bibitem{MVS1}M.V.Safonov; \textit{Boundary estimates for positive solutions to second order elliptic equations},  arXiv:0810.0522v2 [math.AP].
\bibitem{MVS2}M.V.Safonov; \textit{Non-divergence Elliptic Equations
    of Second Order with Unbounded Drift}, Nonlinear partial differential equations and related topics, 211--232, Amer. Math. Soc. Transl. Ser. 2, \textbf{229}, Amer. Math. Soc., Providence, RI, (2010).
\bibitem{WANG}L.Wang; \textit{On the regularity theory of fully nonlinear parabolic equations. II}, Comm. Pure. Appl. Math., \textbf{45} (1992), 141--178.
\end{thebibliography}
\end{document}